\documentclass[12pt]{article}
\usepackage{amsmath,amsfonts,amssymb}
\usepackage{graphicx,color}
\usepackage{theorem}

\theorembodyfont{\upshape}

\newcommand{\ZZ}{{\mathbb Z}}

\newcommand{\CC}{{\mathbb C}}

\def\gg{{\mathfrak{g}}}

\def\hh{{\mathfrak{h}}}

\numberwithin{equation}{section}

\title{Finite dimensional modules for multiloop superalgebra of type  $A(m,n)$ and $C(m)$.}
\author{S.Eswara Rao\\
School of Mathematics\\
Tata Institute of Fundamental Research,\\
Mumbai, India.\\
email: senapati@math.tifr.res.in}
\date{}

\begin{document}

\maketitle

\begin{abstract}
In this article we construct a new class of finite dimensional irreducible modules 
genaralizing evaluation modules for multiloop superalgebra of type $A(m,n)$ and $C(m)$.
We prove that these are the only such modules. In all other cases it is known that such 
modules are evaluation  modules.
\end{abstract}

\section*{Introduction}
In recent years there is some interest in understading finite dimensional modules for the 
Lie algebra $\gg \otimes A$ where $\gg$ is simple finite dimensional Lie algebra and $A$ is 
commutative associative algebra with unit over complex numbers $\CC$. See $[CFK], [RK], [S], [C]$ and
reference there in. Of particular interest is the classification of irreducible finite 
dimensional modules for $\gg \otimes A$. The classification of these modules can be found implictly 
in $[CP]$ in the case $A=\CC [t,t^{-1}].$ They turn out to be evaluation modules at finitely many points. 
Later a direct and simple proof is given in $[E3].$ Similar classification is given for the case 
$A = \CC [t_1^{\pm 1},\cdots, t_d^{\pm 1}]$ the Laurent polynomial algebra in $d$ commuting variables in 
$[E1]$ and $[E2]$. The most genaral case is done in $[CFK]$ with only assumption that $A$ is finitely 
genarated.

In this article we consider $\gg$ to be basic classical Lie superalgebra and $A = \CC [t_1^{\pm 1},\cdots, t_d^{\pm 1}]$. 
The Lie superalgebra $\gg \otimes A$ is called multiloop superalgebra. The main purpose of this article is to classify 
finite dimensional irreducible modules for multiloop superalgebras not covered in $[EZ]$. The classification of 
finite dimensional irreducible modules for multiloop superalgebra is obtained in $[EZ]$ execept for the types 
$A(m,n)$ and $C(m)$. They are all evaluation modules. Here we need to note in the type $A(n,n)$, we have two Lie superalgebras. 
One is the basic classical Lie superalgebra of type $A(n,n)$ which is centerless and the other is $sl(n,n)$ which is 
one dimensional central extension of the earlier algebra. When we say multiloop superalgebra, we include $sl(n, n)$.

As remarked in $[EZ]$ the proof given in $[EZ]$ work for the basic classical Lie superalgebra of type $A(n,n)$. 
Thus we are left with two cases where $\gg = sl(m,n)$ or of type $C(m).$ We note that in these two cases the even part of 
$\gg$ is not semisimple but reductive with one dimensional center. In all other cases the even part is semisimple. 
This is the reason why the proof given in $[EZ]$ does not work for the two cases we consider in this article.\\
The article is organised as follows.

In section 1 we define basic classical Lie superalgebra and recall their classification due to $Kac [K1]$.
We also define multiloop superalgebra. In Section 2 we construct evaluation modules for any multiloop superalgebra. 
It is proved in $[EZ]$ that any finite dimensional irreducible module for multiloop superalgebra is an evaluation 
module execept in the two special cases given above. In Section 3 we describe the Lie superalgebra $sl(m,n)$ and $C(m)$ 
in more detail and describe their root systems. In Section 4 we construct a new class of irreducible finite dimensional 
modules for the two cases above which are not evaluation modules.

This construcction, which is more genaral than evaluation modules, uses the special property that the even part has a 
one dimensional center.

In the last section we prove that any finite dimensional irreducible module for the multiloop superalgebra of the 
$sl(m,n)$ and $C(m)$ is one of the module constructed in Section 4.

\section{Multiloop algebras}
All our vector spaces, algebras and tensor products are over the field of complex numbers $\CC$. $\ZZ_2$ denotes the 
field of two elements $\{{\overline{0},\overline{1}}\}$.

A Lie superalgebra is a $\ZZ_{2^{-}}$ graded vector space $\gg = \gg_{\overline{0}} \oplus \gg_{\overline{1}}$ 
equipped with $\CC-$bilinear form $[,] : \gg \times \gg \rightarrow \gg$, called the Lie super bracket, satisfying the following conditions. 
\begin{enumerate}
\item  $[\gg_{\overline{i}}, \gg_{\overline{j}}] \ {\underline{\subset}} \  \gg_{\overline{i+j}}$

\item $[X,Y] =  - (-1)^{ij} \ [Y,X]$

\item $\big[[X,Y],Z\big] = \big[X,[Y,Z]\big] - (-1)^{ij} \ \big[Y,[X,Z]\big]$
\end{enumerate}
for all homogeneous elements $X \in \gg_{\overline{i}}, Y \in \gg_{\overline{j}}$ and $Z \in \gg_{\overline{k}}$

The subspace $\gg_{\overline{0}}$ is called even and the sbuspace $\gg_{\overline{1}}$ is called odd. It is easy to see that 
$\gg_{\overline{0}}$ is the usual Lie algebra and $\gg_{\overline{1}}$ is $\gg_{\overline{0}}$ module. Suppose $X$ is a 
homogeneous element belonging to $\gg_{\overline{i}}$, then we denote $|X|=i.$

A Lie superalgebra is called basic classical if the algebra is simple, finite dimensional, the even part is reductive 
and carries an even, non-degenerate super symmetric invariant bilinear form. They have been classified by $Kac [K1]$. The 
following is the list of basic classical Lie superalgebra which are not Lie algebras and the decomposion of even part.
$$
\begin{array}{llll}
A(m,n)& A_{m} + A_{n} + \CC& m \geq 0,n \geq 0, m+n \geq 1, m\neq n\\
A(n,n)& A_{n} + A_{n}& n \geq 1\\
B(m,n)& B_{m} + C_{n}& m \geq 0,n \geq 1\\
C(n)& C_{n-1} + \CC& n \geq 3\\
D(m,n)& D_{m} + C_{n}& m \geq 2,n \geq 1\\
D(2,1,a)& D_{2} + A_{1}& a\neq 0, -1\\
F(4)& B_{3}+A_{1}\\
G(3)& G_{2}+A_{1}\\
\end{array}
$$
We will now define multiloop superalgebras. We fix a positive integer $d$ and let $A = \CC [t_1^{\pm 1},\cdots, t_d^{\pm 1}]$ 
be a Laurent polynomial ring in $d$ commuting variables. Let $\gg$ be a basic classical Lie superalgebra. In the case 
$\gg$ is of type $A(n,n), n \geq 1$ we will allow $\gg$ to be $sl(n+1,n+1), n\geq1$. Note that $sl(n+1,n+1)$ has a 
one dimensional center and modulo the center it is basic classical. $\gg \otimes A$ has a natural Lie superalgebra structure.

\paragraph{(1.1) Definition} $\gg \otimes A$ is called multiloop superalgebra.

Note that the universal central extension of $\gg \otimes A$ for $d=1$ are called affine superalgebras. They are widely studied. 
See $[EF]$ and references there in .

\paragraph{(1.2)} The main purpose of this paper is to classify irreducible finite dimensional modules for multiloop 
superalgebras in the following cases.
\begin{enumerate}
\item $\gg = sl(m+1,n+1), m,n\geq 0, m+n \geq1$ 

\item $\gg$ is of ype $C (m), m \geq3$
\end{enumerate}
We will describe these two algebras in more detail in Section 3.

In all other cases the classification has been obtained in $[EZ]$ (Propositions 5.2 and 5.4). The case where $\gg$ is basic 
classical Lie superalgebra of type $A(n,n)$ is similar as noted in Remark 6.7(3) and Remark (7.5) of $[EZ]$.

\section{Evaluation modules}

In this section we will define the all important evaluation modules for multiloop superalgebras. Before that we need to 
recall the notion of highest weight modules for $\gg$ as well as $\gg \otimes A$.

We fix a basic classical Lie superalgebra including the Lie superalgebra $sl(n+1, n+1)$. The following is well known. 
See $[K2]$ for detail.

Let $\gg = N^{-} \oplus \hh \oplus N^{+}$  be the standard decomposition where $\hh$ is a Cartan subalgebra which is 
actully a Lie algebra. For each $\lambda \in \hh^{*}$ there exists an irreducible  $\gg$-module denoted by $V(\lambda)$ which is called 
highest weight module with respect to the above decomposition.

Let $\gg \otimes A = N^{-} \otimes A \oplus \hh \otimes A \oplus N^{+} \otimes A$ be the corresponding decomposition for $\gg \otimes A$.

\paragraph{(2.1)} A module $V$ of $\gg \otimes A$ is called highest weight module if there exists a weight vector $v$ 
(with respect to $\hh$) in $V$ such that 

\begin{enumerate}
\item  $U (\gg \otimes A) v = V$

\item $N^{+} \otimes A v = 0$

\item $U (\hh \otimes A)v = \CC v.$

\end{enumerate}

Here $U$ denote the universal enveloping algebra. By standard arguments we can show that for each $\psi \in (\hh \otimes A)^{*}$ 
there exists a unique irreducible highest weight module denoted by $W(\psi)$.

Now we recall the evaluation modules. For each $i, 1\leq i\leq d,$ let $N_{i}$ be a positive integer. Let 
$\underline{a}_{i}= (a_{i1}, a_{i2}, \cdots a_{{iN}_{i}})$ be non-zero distinct complex numbers.\\
Let $N=N_{1} N_{2} \cdots N_{d}.$ Let $I= (i_{1}, \cdots, i_{d})$ where $1\leq i_{j} \leq N_{j}$\\
There are $N$ of them. Let $I_{1}, \cdots I_{N}$ be some order.\\ Let $\underline{m}=(m_{1}, \cdots m_{d}) \in \ZZ^{d}$
and define $a^{\underline{m}}_{I} = a^{m_1}_{1i_1} \cdots a^{m_d}_{di_d}$.\\
Let $t^{\underline{m}}= t^{m_1}_{1} \cdots t^{m_d}_{d} \in A$.

Let $G$ be any Lie superalgebra and let $\varphi$ be a algebra homomorphism defined by
$$
\begin{array}{llll}
\varphi :G \otimes A \rightarrow \oplus G &=& G_{N} (N copies)\\
\varphi(X \otimes t^{\underline{m}}) &=& (a^{\underline{m}}_{I_{1}} X, \cdots a^{\underline{m}}_{I_{N}} X)
\end{array}
$$
Where $X \in G$ and $\underline{m} \in \ZZ^{d}$ \\
Let $P_{j} (t_{j})= \displaystyle{\prod^{N_{j}}_{k=1}} (t_{j}-a_{jk})$ and note that $P_{j}$ is a polynomial in $t_{j}$ 
with nonzero distinct roots.\\
Let $I$ be the co-finite ideal of $A$ genarated by $P_{1}(t_{1}), \cdots P_{d}(t_{d}).$

\paragraph{Lemma (2.2)}
\begin{enumerate}
 \item $\varphi$ is surjective

 \item Ker $\varphi = G \otimes I$ and $G \otimes A/I \cong G_{N}.$
\end{enumerate}
\paragraph{Proof} The proof is similar to the Lie algebra case and can be found in Lemma 3.11 of $[E2]$. A more 
transparent proof can be found in Proposition 2.2 of $[B]$.

Let $V(\lambda_{1}), \cdots V(\lambda_{N})$ be irreducible highest weight modules for $\gg$. Then 
$V (\underline{\lambda}, \underline{a}) = \displaystyle{\bigotimes^{N}_{i=1}} V(\lambda_{i})$ is an irreducible highest weight 
module for the multiloop algebra $\gg \otimes A$ via the map $\varphi$ taking $G = \gg $.\\
Let $\psi : \hh \otimes A \rightarrow \CC$ be defined by \\
$\psi (h \otimes t^{\underline{m}}) = \displaystyle{\sum^{N}_{j=1}} a^{\underline{m}}_{I_j} \lambda_{j} (h), h \in \hh$\\
Then it is easy to see that $W(\psi) \cong V({\underline{\lambda}},{\underline{a}}).$

\paragraph{(2.3) Definition} The modules $V({\underline{\lambda}},{\underline{a}})$ are called evaluation modules.\\

\noindent
Note that $V({\underline{\lambda}},{\underline{a}})$ is finite dimensional if and only if each $V(\lambda_{i})$ is finite 
dimensional. The condition for $V(\lambda_{i})$ to be finite dimensional is given in $[K2]$.

\paragraph{(2.4) Proposition} Suppose $\gg$ is not of type $A(m,n)$ or $C(m).$ Then any finite dimensional irreducible 
module is an evaluation module.

\paragraph{Proof} See Propositions 5.2 and 5.4 of $[EZ].$

In the case $\gg$ is of type $A(m,n)$ or $C(m)$, there exists irreducilble finite dimensional modules which are not evaluation 
modules which we will see in Section 4. Here we leave out that case where $\gg$ is basic classical Lie superalgebra of 
type $A(n+1,n+1).$ In this case the classification is similar to the one given in Proposition 2.4.

\section{Lie superalgebra $sl(m+1, n+1)$ and $C(m).$}
\quad In this section we assume $\gg$ is $sl(m+1, n+1), (m\geq 0, n \geq 0, m+n \geq 1)$ or $\gg$ is of type $C(m), m \geq 3.$ 
We will now describe these two cases of Lie superalgebras and their root systems. We follow $Kac [K2]$.

\paragraph{Case 1} $\gg = sl(m, n)$\\
Let $V=V_{\overline{0}} \oplus V_{\overline{1}}$ be $\ZZ_{2}$- graded vector space and let $\dim V_{\overline{0}} = m$ and 
$\dim V_{\overline{1}} = n$. The $\ZZ_{2}$- gradation on $V$ naturally induces a $\ZZ_{2}$- gradation on \\
$End(V) = (End (V))_{\overline{0}} \oplus (End (V))_{\overline{1}}$\\
by letting 
$End(V)_{j} = \{f \in End(V) : f (V_{k})_{\subseteq} V_{k+j} for \ all \ k \in \ZZ_{2}\}$\\
Then $End (V)$ becomes a Lie superalgebra by defineing Lie super bracket\\
$[f,g]=f \circ g - (-1)^{ij} g \circ f$ for all $f \in (End(V))_{i}$ and $g\in (End(V))_{j}$.\\
We let $gl(m,n) = End(V).$ By fixing a basis for $V_{\overline{0}}$ and $V_{\overline{1}}$ we can write $f \in End(V)$ as

$$
f=
\begin{pmatrix}
A&B\\
C&D
\end{pmatrix}
$$
Where $A$ is $m \times m$ matrix, $B$ is $n \times m$ matrix, $C$ is $m \times n$ matrix and $D$ is $n \times n$ matrix. 
It is easy to see that 
$$
\begin{pmatrix}
A&0\\
0&D
\end{pmatrix} \in \big(End (V)\big)_{\overline{0}}
$$
and
$$
\begin{pmatrix}
0&B\\
C&0
\end{pmatrix} \in \big(End (V)\big)_{\overline{1}}
$$
We now define super trace of $f$ denoted by strf= trace A -trace D.\\
Define 
$sl(m,n) = \{X \in gl(m,n)|strX=0\}$\\
It is known that $sl(m,n)$ is simple Lie superalgebra if $m\neq n$. \\When $m=n, sl(n,n)$ has a  one dimensional center if 
$n\geq 2.$\\
Define
$$
\begin{array}{llllll}
\gg_{0}= \left\{ \begin{pmatrix}A&0\\0&D\end{pmatrix}; A \mbox{ is} \ m \times m  \ \mbox{matrix  and  D  is} 
\ n \times n \ \mbox{ matrix} \right\} \\ [4mm]                                     
\gg_{+1}= \left\{ \begin{pmatrix}0&B\\0&0\end{pmatrix}; B \mbox{ is} \  n \times m   \ \mbox{matrix} \right\}\\
\gg_{+1}= \left\{ \begin{pmatrix}0&0\\C&0\end{pmatrix}; C \mbox{ is} \  m \times n   \ \mbox{matrix} \right\}                                        
\end{array}
$$
Then $\gg_{\overline{0}}= \gg_{0}, \gg_{\overline{1}}= \gg_{-1} \oplus \gg_{+1}$\\
It is easy to check that\\
$
(3.1) \ [\gg_{0},\gg_{-1}] \subseteq \gg_{-1},
[\gg_{0},\gg_{+1}] \subseteq \gg_{+1},
[\gg_{-1},\gg_{+1}] \subseteq \gg_{0},
[\gg_{\pm 1},\gg_{\pm 1}] = 0.
$

Thus $\gg = \gg_{-1} \oplus \gg_{0} \oplus \gg_{+1}$ is a $\ZZ$-grading with $\gg_{n} = 0$ for $|n| >1.$\\
We can further check that
$
[\gg_{{\overline{1}}} , \gg_{{\overline{1}}}] \nsubseteq  [\gg_{{\overline{0}}} , \gg_{{\overline{0}}}]\\
$
This is the reason why the proof in $[EZ]$ does not work for type $A(m,n)$ and $C(m).$

In all other cases (except in the case 1.2(1) and 1.2(2)) $\gg_{\overline{0}}$ is semisimple and hence 
$
[\gg_{{\overline{1}}} , \gg_{{\overline{1}}}] \subseteq \gg_{\overline{0}}= [\gg_{{\overline{0}}} , \gg_{{\overline{0}}}].\\
$
Let $\hh$ be a Cartan subalgebra of $\gg$ which is contained in the even part $\gg_{\overline{0}.}$ We can take $\hh$ to be 
the diognal matricas in $\gg_{\overline{0}.}$ Let $\epsilon_{i}$ be the ith projection from $\hh$ to $\CC$. 
Let $\delta_{i} = \epsilon_{m+i}$ for $1 \leq i \leq n.$ Then the root system of $sl(m,n)$ can be describe in the following way
$$
\begin{array}{llll}
\varDelta_{0} &=& \left\{ \epsilon_{i} - \epsilon_{j}, \delta_{i} - \delta_{j}, i \neq j \right\}\\
\varDelta_{1} &=& \left\{ \pm (\epsilon_{i} - \delta_{j})\right\}
\end{array}
$$
Here $\varDelta_{0}$ denotes the even roots and $\varDelta_{1}$ denotes the odd roots.\\
Let $\varDelta_{1}^{+}= \{ \epsilon_{i} - \delta_{j} \}$ and $\varDelta_{1}^{-}= \{ -(\epsilon_{i} - \delta_{j}) \}$\\
Then it is easy to see that the roots of the space $\gg_{+1}$ are $\varDelta_{1}^{+}$ similarly for  
$\gg_{-1}$ are $\varDelta_{1}^{-}$. From the above description of the root system one can easily see that (3.1) holds.
\paragraph{Case 2} $\gg = C(m), m \geq 3.$

We will now describe the Lie superalgebra $C(m)$ which is known to be a subalgebra of $sl(2,2m-2)$. We denote $A^T$ for the
transpose of the matrix A.\\
Let 
$$
\gg_0 =\left\{\begin{array}{cp{2.5in}}
\left(
\begin{array}{lcc|ccl}
  \alpha & 0 &&0&0 &\\
0 &  -\alpha &&0 &0 \\ \hline
0&0&& A& B \\
0&0&& C&-A^T \\
\end{array}
\right) : &  
$\alpha \in \CC, A, B, C$ \ are \ $(m-1)\times (m-1)$
\ matrices, $B^T=B, C^T =C$ \end{array}
\right\}
$$

$$
\gg_{+1} =\left\{\left(\\
\begin{tabular}{lcc|ccl}
0 & 0 &&$A$&$B$ &\\
0 &  0&&0 &0 \\ \hline
0&$-B^{T}$&& 0& 0 \\
0&A&& 0&0 \\
\end{tabular}
\right) ; A \mbox{ and} \ B \ \mbox{ are} \ (m-1)\times 1
\ \mbox{matrices} \right\}
$$

$$
\gg_{-1} =\left\{\left(\\
\begin{tabular}{lcc|ccl}
0 & 0 &&0&0 &\\
0 &  0&&$A$ &$B$ \\ \hline
$-B^{T}$&0&& 0& 0 \\
$A$&0&& 0&0 \\
\end{tabular}
\right) ; A, B \ \mbox{ are} \ (m-1)\times 1
\ \mbox{matrices} \right\}
$$

Then it is easy to check \\
$
(3.2)\  [\gg_{0},\gg_{+1}] \subseteq \gg_{+1},
[\gg_{0},\gg_{-1}] \subseteq \gg_{-1},
[\gg_{-1},\gg_{+1}] \subseteq \gg_{0},
[\gg_{\pm 1},\gg_{\pm 1}] = 0.
$

From this we can see that
$\gg=\gg_{-1} \oplus \gg_{0} + \gg_{1}$ is a $\ZZ$-gradation with $\gg_{n}=0$ for $|n|>1$.\\
Further we have 
$$
\gg_{0}=\gg_{\overline{0}}, \gg_{\overline{1}} =\gg_{-1} \oplus \gg_{+1}
$$
We can also check that 
$$
[\gg_{\overline{1}}, \gg_{\overline{1}} ] \nsubseteq [\gg_{\overline{0}},\gg_{\overline{0}}]
$$
We will now describe the root system of $C(m)$. As in the earlier case we can take the Cartan subalgebra $\hh$ to 
be diognal matricas in $\gg_{\overline{0}}$. Let $\epsilon_{i}$ be the ith projection from $\hh$ to $\CC$.\\
Let $\delta_{1} = \epsilon_{3}, \cdots \delta_{m-1} = \epsilon_{m+1}$\\
Let $\Delta_{0}, \Delta_{1}, \Delta_{1}^{+}, \Delta_{1}^{-}$ be even, odd, positive odd and negative odd roots respectively.\\
Then
$$
\begin{array}{llll}
\Delta_{0}&=& \{\pm 2 \delta_{i}, \pm \delta_{i} \pm \delta_{j}, i \neq j \}\\
\Delta_{1}&=& \{\pm \epsilon_{1} \pm \delta_{i} \}\\
\Delta_{1}^{+}&=& \{\epsilon_{1} \pm \delta_{i} \}\\
\Delta_{1}^{-}&=& \{-\epsilon_{1} \pm \delta_{i} \}\\
\end{array}
$$
It is easy to verify that $\Delta_{1}^{+}$ is the set of roots for $\gg_{+1}$ and $\Delta_{1}^{-}$ is the set of roots of 
$\gg_{-1}$. One can also verify (3.2) useing above description.

\section{A new class of finite dimensional modules for multiloop superalgebra of type $A(m,n)$ and $C(m)$.}
Let $\gg=sl(m+1,n+1), m\geq 0, n\geq0, m+n\geq 1)$ or $\gg=C(m), m\geq3.$ Recall that $\gg_{\overline{0}}$ is reductive but not 
semisimple. Let $\gg_{ss}$ be the semisimple part of  $\gg_{\overline{0}}$ and $\CC z$ be the one dimensional center of  
$\gg_{\overline{0}}$. So that we have 
$$
\gg_{\overline{0}} = \gg_{ss} \oplus \CC z
$$

Let $V$ be finite dimensional evaluation modules for $\gg_{ss} \otimes A.$ Thus there exists a co-finite ideal $I^{'}$ of 
$A$ genarated by polynomial $P_{j}^{'} (t_{j}) = \displaystyle{\prod^{N_{j}}_{k=1}} (t_{j}- a_{jk})$ where $a_{jk}$ and 
$N_{j}$ are given in section 3 such that $\gg_{ss} \otimes I^{'} . V=0$\\
Further $V \cong W(\psi)$ where $\psi \in (\hh \otimes A)^{*}$ and $\psi (h \otimes t^{\underline{m}}) = \sum 
a_{I_{j}}^{m} \lambda_{j} (h).$ (See section 3 for details).\\
We will allow $\psi$ to be zero and note that $W(0)$ is trivial one dimensional module.\\
Let $b_{jk}$ be some positive integers for $1\leq j \leq d$ and $1 \leq k \leq N_{j}$. Let 
 $P_{j}(t_{j}) = \displaystyle{\prod^{N_{j}}_{k=1}} (t_{j}- a_{jk})^{b_{jk}}.$ Let $I$ be the ideal genarated by 
$P_{j} (t_{j}), j=1,2,\cdots d.$ clearly $I \subseteq I^{'}.$ Let $\lambda \in (\CC z \otimes A)^{*}$ such  that 
$\lambda (z \otimes I)=0$. Now consider $W(\psi)$ as $D_{\circ}: = \gg_{\overline{0}} \otimes A/I = \gg_{ss} \otimes A/I \oplus \CC z \otimes A/I$
module. Let $D_{+} = \gg_{+1} \otimes A/I (D_{-} = \gg_{-1} \otimes A/I) $ act on $W(\psi)$ trivially. Since 
$[D_{\circ}, D_{+}] \subseteq D_{+}$ and $[D_{+}, D_{+}] = 0$, we see that $W(\psi)$ is a well defined module for 
$D_{\circ} \oplus D_{+}$. Now consider the induced module for $\gg \otimes A/I$
$$
M(\psi, \lambda) = U(\gg \otimes A/I) \displaystyle{\bigotimes_{D_{\circ} + D_{+}}} W(\psi)
$$
and as vector space $\cong U(D_{-}) \otimes W(\psi)$. Since $D_{-}$ is finite dimensional and odd, by the PBW theorem for 
Lie superalgebras, we conclude that $M(\psi, \lambda)$ is finite dimensional. By standard argument we see $M(\psi, \lambda)$ 
has a unique irreducible quotient say $V(\psi, \lambda)$. By the surjective map $\gg \otimes A \rightarrow \gg \otimes A/I$ 
we see that  $V(\psi, \lambda)$ is an irreducible finite dimensional module for $\gg \otimes A.$ It is not hard to see that 
$V(\psi, \lambda)$  is not an evaluation module if $\lambda(z \otimes I^{'}) \neq 0.$ In the case when 
$\lambda(z \otimes I^{'})= 0$, $V(\psi, \lambda)$ is an evaluation module. Follows from Lemma 5.1.

\section{Classification Theorem.}
In this section we recall some genaralities on highest weight modules for multiloop superalgebras. We will also state and 
prove our main result on the classification. 

Let $\gg$ be any finite dimensional Lie superalgebra as mentioned in \break Section 1.

Consider the standard decomposition.

$\gg \otimes A= N^{-} \otimes A \oplus \hh \otimes A \oplus N^{+} \otimes A$ \\
Note that $\hh \otimes A$ is a Lie algebra (abelian). Let $\psi \in (\hh \otimes A)^{*}$ and $\CC_{\psi}$ be one dimensional 
representation of $\hh \otimes A$ via $\psi$. Let $N^{+} \otimes A$ act trivially on $\CC_{\psi}$. Consider the Verma module.
$$
M(\psi)= U(\gg \otimes A) \bigotimes_{N^{+} \otimes A + \hh \otimes A} \CC_{\psi}
$$
By standard arguments we see that $M(\psi)$ has a unique irreducible quotient say $V (\psi)$. $V (\psi)$ need not have 
finite dimensional weight spaces (with respect to $\hh$).

\paragraph{Lemma 5.1}
$V (\psi)$ has finite dimensional weight spaces if and only if $\psi$ factors through $\hh \otimes A/I$ for some co-finite 
ideal $I$ of $A$. In this case $\gg \otimes I. V(\psi)=0$.

\paragraph{Proof} Follows from Lemma 4.7 and Remark 4.8 of $[EZ]$. See Lemma 3.7 of $[EZ]$ for the original proof. Note that 
$\gg^{'}$ that occurs in $[EZ]$ is $[\gg,\gg]$ and equals to $\gg$ in our case. 

\paragraph{Remark 5.2} From the proof of Lemma 3.7 of $[E2]$ it follows that the co-finite $I$ can be chosen to to be 
genarated by polynomials $P_{1}, \cdots P_{d}$ in variables $t_{1}, \cdots t_{d}$ respectively. We can also assume that they have 
non-zero roots.

From now onwards we assume $\gg = sl(m+1, n+1)$ or $C(m).$

\paragraph{Theorem 5.3} Any finite dimensional irreducible module for $\gg \otimes A$ is isomorphic to $V(\psi, \lambda)$ as 
defined in Section 3.

\paragraph{Proof} Let $V$ be finite dimensional irreducible module for $\gg \otimes A$. Since $V$ is finite dimensional 
$V \cong V(\phi)$ for some $\phi \in (\hh \otimes A)^{*}$.\\
Let $v$ be the highest weight vector. Recall that $\gg_{ss}$ is the semisimple part of $\gg_{\overline{0}}$. Let 
$M=U(\gg_{ss} \otimes A)v$ which is finite dimensional $\gg_{ss} \otimes A$ module. In fact it is $\gg_{\overline{0}} \otimes A$ 
module as the additional center $z \otimes A$ act as scalars.

\paragraph{Claim} $M$ is irreducible as $\gg_{ss} \otimes A$ module.

To see the claim, recall that $\gg = \gg_{-1} \oplus \gg_{\overline{0}} \oplus \gg_{+1}$ and $\gg_{+1}$ is sum of positive 
odd root spaces. Let $\gg_{ss}=\gg_{ss}^{-} \oplus \hh_{ss} \otimes \gg_{ss}^{+}$ be the standard triangular decompontion. Since
$v$ is highest weight vector we have 
$$
M=U(\gg_{ss}^{-} \otimes A)v
$$
Let $w \in M$ be a weight vector of weight $\mu_{1}$ and let $\phi | \hh = \mu.$ Then $w= X v$ for 
$X \in U(\gg_{ss}^{-} \otimes A)$\\
No loss of genarality by assuming  $X$ is a monomial and of weight $-\beta$ where $\beta$ is non-negative linear combination of 
simple roots comeing from $\gg_{ss}$.\\
We have $\mu_{1}=\mu-\beta$.\\
Since $V(\phi)$ is $\gg \otimes A$-irreducible there exists $Y \in U \big((\gg_{+1} \oplus \gg^{+}_{ss})\otimes A \big)$
such that $Yw = v.$ We can assume $Y$ is monomial and of weight $\beta$. By looking at the root systems, this is not possible 
unless $Y \in U (\gg^{+}_{ss} \otimes A).$ This proves that $M$ is $\gg_{ss} \otimes A$ irreducible and hence the claim.\\
Now by Lemma 4.1 and Remark 4.2 there exists a co-finite ideal $I$ of $A$ such that 
$$
\gg \otimes I. V(\phi) =0.
$$ 
We further can choose polynomials $P_{1}, \cdots P_{d}$ in variables $t_{1}, \cdots t_{d}$ such that $I$ is genarated by 
$P_{1}, \cdots P_{d}$.\\
Let $P_{j}(t_{j}) = \displaystyle{\prod^{N_{j}}_{j=1}} (t_{j} - a_{jk})^{b_{jk}}$ where for each 
$j, a_{j1}, \cdots a_{jN_{j}}$ are distinct non-zero complex numbers. $b_{jk}, N_{j}$ are some positive  integers.\\
Let $P_{j}^{'} (t_{j}) = \displaystyle{\prod^{N_{j}}_{j=1}} (t_{j} - a_{jk})$ and let $I^{'}$ be the ideal genarated by 
$P_{1}^{'}, \cdots P_{d}^{'}$ clearly $I \subseteq I^{'}.$ Now from the classification of finite dimensional irreducible modules 
for $\gg_{ss} \otimes A$ (See Remark 5.5 of $[EZ]$. There the proof is given for simple Lie-algebras but the same proof works
for any finite dimensional semisimple Lie algebra) it follows that 
$$
(\gg_{ss} \otimes I^{'}) M=0
$$ 
The space $\CC z \otimes A$, which is central in $\gg_{\overline{0}} \otimes A$, act on $M$ as scalars via the map $\phi$ and 
factors through $\CC z \otimes A/I$. Further $\gg_{+1} \otimes A$ acts trivially on $M$ as 
$[\gg_{+1} \otimes A, \gg_{ss} \otimes A] \subseteq \gg_{+1} \otimes A$\\
Thus $M$ satisfy all the properties of $W(\psi)$ in section 3. Thus by the universal property of induced modules, 
$V(\phi)$ is the quotient of $U(\gg \otimes A) \displaystyle{\bigotimes_{D_{\circ} + D_{1}}} M.$\\
Hence $V(\phi) \cong V(\psi, \lambda)$ for some $\psi$ and $\lambda$ as in Section 3.\\
This completes the proof of the Theorem.

\end{document}